\RequirePackage{ifpdf}
\ifpdf 
\documentclass[pdftex]{sigma}
\else
\documentclass{sigma}
\fi

\def\a{{\alpha}}
\def\b{{\beta}}
\def\o{\over}
\def\P{\Psi}
\def\lm{\lambda}
\def\d{\delta}
\def\g{\gamma}
\def\D{{\Delta}}

\numberwithin{equation}{section}

\begin{document}

\allowdisplaybreaks

\renewcommand{\thefootnote}{$\star$}

\renewcommand{\PaperNumber}{089}

\FirstPageHeading

\ShortArticleName{A Probablistic Origin for a New Class of Bivariate Polynomials}
\ArticleName{A Probablistic Origin for a New Class\\ of Bivariate Polynomials\footnote{This paper is a contribution to the Special
Issue on Dunkl Operators and Related Topics. The full collection
is available at
\href{http://www.emis.de/journals/SIGMA/Dunkl_operators.html}{http://www.emis.de/journals/SIGMA/Dunkl\_{}operators.html}}}

\renewcommand{\thefootnote}{\arabic{footnote}}
\setcounter{footnote}{0}

\Author{Michael R. HOARE and Mizan RAHMAN\,\footnote{Supported partially by an NSERC Grant \#A6197.}}

\AuthorNameForHeading{M.R. Hoare and M. Rahman}

\Address{School of Mathematics and Statistics, Carleton University, Ottawa, ON K1S 5B6, Canada}
\Email{\href{mailto:mrahman@math.carleton.ca}{mrahman@math.carleton.ca}}

\ArticleDates{Received September 15, 2008, in f\/inal form December 15,
2008; Published online December 19, 2008}

\Abstract{We present here a probabilistic approach to the generation of new polyno\-mials in two discrete variables. This extends our earlier work on the `classical' orthogonal polynomials in a previously unexplored direction, resulting in the discovery of an exactly soluble eigenvalue problem corresponding to a bivariate Markov chain with a transition kernel formed by a convolution of simple binomial and trinomial distributions. The solution of the relevant eigenfunction problem, giving the spectral resolution of the kernel, leads to what we believe to be a new class of orthogonal polynomials in two discrete variables. Possibilities for the extension of this approach are discussed.}


\Classification{33C45; 60J05}

\section{Introduction}\label{section1}

Some thirty years ago we published several papers \cite{6, 10,11,12, 15, 16} in which we described a~class of statistical models which gave rise to the `classical' orthogonal polynomials and a variety of associated formulas which had previously been known only in the abstract. The key to this was to def\/ine simple Markov chains using certain `Urn models' and variants of Bernoulli trials, whose transition kernels provided soluble eigenvalue problems that in turn yielded the polynomials of interest as eigenfunctions. This was carried out for both continuous and discrete variables and led to a scheme in which the quintet of discrete-single-variable orthogonal polynomials (Hahn, Gonin, Krawtchouk, Meixner, Charlier) could be inter-related by suitable limits and substitutions. These in turn underlay the better-known continuous-variable sets (Laguerre, Jacobi). As well as solving the def\/ining eigenvalue problems, we were able to discover results relating, amongst others, to `ladder-operators', the Factorization Method, and dual polynomials, some ef\/fectively new, others known in disguised form in the `Bateman-project' era.

The success of these methods opens up a more interesting prospect, namely that of actually discovering {\it new} polynomial systems through statistical models, or perhaps, less ambitiously, distinguishing specially important cases within the generality of those already known. While the scope for this in the case of the single variable was predictably limited, a whole new opportunity presents itself when functions of more than one variable are considered. While a plethora of special functions of several variables has emerged in the last decades, in various stages of genera\-li\-ty, few have been related to tangible structures, such as might have their origin in physical, statistical or combinatorial models. With the experience of the single variable results, there is reason to hope that the eigenvalue problems based on one of these might lead to a~specially interesting class of solutions illuminating the hitherto rather unstructured world of multivariate special functions.

In this paper we shall describe our f\/irst steps in this direction. We consider an extension to two discrete variables of the `Cumulative Bernoulli Trials', f\/irst described in 1983 \cite{10}, and show that this leads to a soluble eigenvalue problem via a transition-kernel involving the bi- and trinomial distributions. The resulting eigenfunctions prove to be an adaptation of the `$9-j$ symbols' known in theoretical physics, where they form a central idea in the theory of angular momentum. Once again a mathematical structure proves to pervade the natural world in the most surprising places. Here we shall concentrate on the mathematical content of this result, exploring the probability of a non-trivial explicitly soluble multivariate Markov Chain, possibly for the f\/irst time, is a most satisfying result in itself. See \cite{5, 11,12} for related ideas.

{\bf The statistical model.}
The idea behind `Cumulative Bernoulli Trials' (CBTs) is extraordinarily simple, for all that it seems to have appeared only as late as the 1970s. Whereas ordinary Bernoulli Trials (BTs) represent the outcome of a set of success/failure events with an assigned success probability, Cumulative Bernoulli Trials (CBTs) allow for the possibility that successes in an initial trial can be `saved' while further trials are carried out with the `failures' to increase the number of `successes'. The elementary properties of such a scheme are detailed in~\cite{10}. The simplest realization of the CBT process for descriptive purposes involves trials by `throwing' a~set of dice with a def\/ined success criterion.

Consider the case of six `Poker dice' with faces marked as usual from `Ace' to `nine'. Then the occurrence of $k$ `Ace' on a single throw has the binomial probability
\begin{gather}
b(k,N;\a) = {N\choose k} \a^k (1-\a)^{N-k},\label{1.1}
\end{gather}
where in the example given $N = 6$ and $\a = 1/6$. In \cite{14} we showed that, if the `failed' dice are re-thrown $n$ times with the successes saved, the probability of $i$ successes altogether is: $b^{(n)}(i,N;\a) = b(i, N;1-(1-\a)^n)$, a result that could be generalized to the case where $\a$ is not constant on successive throws. This result is not needed in the present paper, but will serve to motivate the present results.

{\bf Realization of a bivariate Markov chain.}
For convenience we shall continue in the language of Bernoulli Trials with dice, though this is not essential to the structure of the problem. Consider thus a set of $N$ dice, each with a given number of faces, $n$. The faces can be marked in any way, but of these two (or possibly more), for example {\it red} and {\it black} colours, are designated `interesting' outcomes and are scored as `successes' when all or some of $N$ dice are thrown. The chain can be described. Note that they are anti-correlated, since obviously getting more {\it red} reduces the chances of getting many {\it black} and vice versa. Thus the Markov chain is a non-trivial extension from the single variable, not simply a multi-dimensional case. The sides that turn up which are not {\it red} or {\it black} are `failures' and can be called `blanks'. The probabilities of getting {\it red} or {\it black} on throwing a single die are assigned as $\a_1$, $\a_2$ and in the case of the dice are related to $n$ in an obvious way.

{\bf Example.} Standard `Poker dice' $n=6$; {\it black} = `Ace', {\it red} = `King', $\a_1 = \a_2 = 1/6$, $N = 5$. Possible ranges $0 \le i_1 \le N - i_2$; $0 \le i_2 \le N-i_1$. Thus the possible {\it states} of the system are: (0,0), (0,1), (0,2), (0,3), (0,4), (0,5), (1,0), (1,1), (1,2), (1,3), (1,4), (2,0), (2,1), (2,2), (2,3), (3,0), (3,1), (3,2), (4,0), (4,1), (5,0), 21 in all, constituting the {\it state space}.

Returning now to a general $N, \a_1, \a_2$, consider an initial $(i_1, i_2)$ with $N-i_1-i_2$ `blanks'.

{\bf Step 1.}  The $i_1$ `black' dice are thrown, giving $k_1 \le i_1$ `black successes', the $i_2$ `red' dice are thrown, giving $k_2 \le i_2$ `red successes', and these are saved. The probabilities are the respective binomials $b(k_1, i_1; \a_1)$ and $b(k_2, i_2; \a_2)$.

{\bf Step 2.}  The $i_1 + i_2-k_1 - k_2$ `blanks' from the previous step are added to the $N-i_1 - i_2$ original `blanks' giving $N-k_1 -k_2$ `blanks' in all.

{\bf Step 3.}  The collected $N-k_1-k_2$ `blanks' are thrown and the `red' and `black' successes recorded as $p_1$, $p_2$ respectively. These can be at dif\/ferent probabilities $\b_1$, $\b_2$. The outcome is given by the trinomial $b_2(p_1, p_2, N-k_1 - k_2; \b_1,\b_2)$. Evidently:
\[
0 \le p_1 \le N-k_1 - k_2 - p_2; \qquad 0 \le p_2 \le N-k_1 - k_2 - p_1.
\]
Now redef\/ine $p_1 = j_1 - k_1$; $p_2 = j_2 - k_2$ where now
\[
k_1 \le j_1 \le N-j_2; \qquad k_2 \le j_2 \le N-j_1.
\]

{\bf Step 4.}  Combine the `successes' $j_1 - k_1$ and $j_2 - k_2$ with the `successes' held over from Step~1. The `score' of `successes' will now be $(j_1, j_2)$, with $N-j_1 - j_2$ `blanks' and the above process will have led to the transition $(i_1, i_2) \to (j_1, j_2)$. The transition probability for this will be the kernel $K(j_1, j_2; i_1, i_2)$, which clearly def\/ines the {\it transition kernel} of a bivariate Markov chain giving the probability of arriving at state $(j_1, j_2)$ from state $(i_1, i_2)$.

{\bf Step 5.}  Repeat the whole process Steps 1, to 4, to generate the chain.

In the `Poker dice' example the transition matrix will have $21^2 = 441$ elements, many of which will, however, be zero. The sequence of states is Markovian by virtue of the `memory' carried over from the original to f\/inal states as guaranteed by the sequence above. The process will be def\/ined by the {\it transition kernel} $K(j_1, j_2; i_1, i_2)$ the form of which follows by summing all possible pathways, leading to the convolution:
\begin{gather}
K(j_1, j_2; i_1, i_2) = \sum^{\min(i_1,j_1)}_{k_1 = 0} \sum^{\min (i_2,j_2)}_{k_2=0} b (k_1, i_1;\a_1)
 \nonumber\\
 \phantom{K(j_1, j_2; i_1, i_2) =}{} \times b(k_2, i_2;\a_2) b_2 (j_1 - k_1, j_2-k_2, N-k_1-k_2;\b_1, \b_2),\label{1.2}
\end{gather}
where $b(\cdot, \cdot;\cdot)$ is the simple binomial as before and the trinomial $b_2(\cdot, \cdot,\cdot;\cdot,\cdot)$ is
\begin{gather}
b_2 (i_1, i_2, N;p,q) = p^{i_1} q^{i_2} [1-p-q]^{N-i_1 - i_2} \left({N!\o i_1!i_2!(N-i_1-i_2)!}\right).
\label{1.3}
\end{gather}
Thus the kernel is explicitly
\begin{gather}
K(j_1, j_2; i_1, i_2) = i_1!i_2!\b^{j_1}_1 \b_2^{j_2} [1-\b_1-\b_2]^{N-j_1-j_2} {(1-\a_1)^{i_1}(1-\a_2)^{i_2}\o
(N-j_1 - j_2)!}\nonumber\\
\phantom{K(j_1, j_2; i_1, i_2) =}{}  \times \sum^{\min(i_1,j_1)}_{k_1=0}  \sum^{\min (i_2,j_2)}_{k_2=0}
 \left({\a_1\o 1-\a_1}\right)^{k_1} \left({\a_2\o 1-\a_2}\right)^{k_2}\nonumber\\
\phantom{K(j_1, j_2; i_1, i_2) =}{}  \times {1\o \b_1^{k_1}\b_2^{k_2}}  {(N-k_1-k_2)!\o
(i_1-k_1)!(i_2-k_2)!(j_1-k_1)!(j_2-k_2)!k_1!k_2!}.\label{1.4}
\end{gather}

\begin{figure}[t]

\centerline{\includegraphics{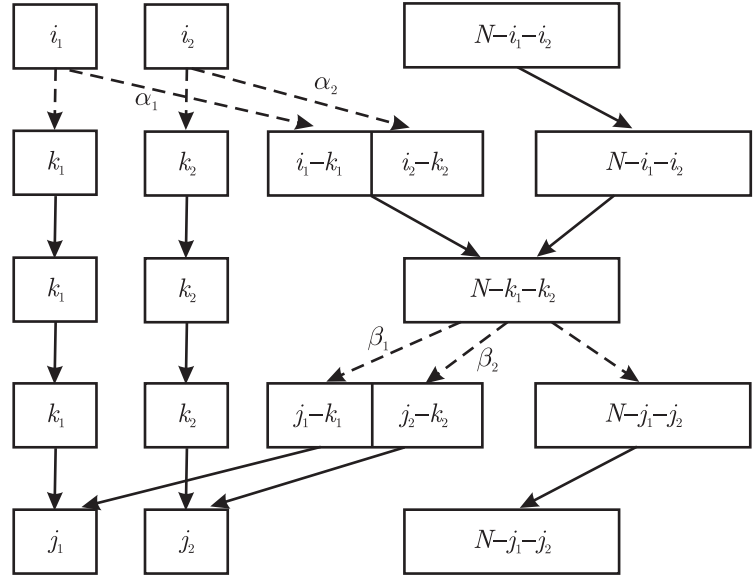}}

\caption{Schematic representation of the Markov chain for cumulative binomial and trinomial trials. Note that the dotted lines refer to stochastic outcomes while the solid arrows indicate counts carried forward.}
\label{fig1}
\end{figure}

Fig.~\ref{fig1} will be helpful in following the steps outlined above. Having established the probabilistic background to the kernel $K$ we shall now consider its eigenvalue problem on the discrete state-space $(i_1, i_2)$ with $0 \le i_i$, $i_2 \le N$.
\begin{gather}
\sum_{j_1}\sum_{j_2} K(i_1, i_2; j_1, j_2) \P_{m,n}(j_1,j_2) = \lm_{m,n}\P_{m,n}(i_1, i_2).\label{1.5}
\end{gather}
This is the form we shall take as origin of the eigenvalue problem to be investigated. For the moment we may note the crucial property that $K$ is a function of four parameters $0 \le \a_1$, $\a_2$, $\b_1$, $\b_2 \le 1$. Before embarking on this, however, we shall need to spend some time on an excursion, placing the problem in the context of present-day theory of multivariate orthogonal polynomials. Only then can we tackle the eigenfunction problem for~$K$.

\section{Multivariate orthogonal polynomials of a discrete variable}\label{section2}

The modern era of single-variable classical orthogonal polynomials probably began with Wilson's~\cite{22,23} observation that Wilson's $6-j$ symbols which are known as the Racah~\cite{18} coef\/f\/icients, are constant multiples of the $_4F_3$ polynomial:
\begin{gather}
R_n(x):= {}_4F_3\left[\begin{array}{c} -n, n+\a + \b+1, -x, x+\g - N\\ \a+1, -N, \b + \g +1\end{array} ;1\right],\label{2.1}
\end{gather}
$x$, $n = 0, 1, \ldots, N$, and that they have a $q$-analogue of the form
\begin{gather}
{}_4\Phi_3\left[\begin{array}{c}
q^{-n}, abq^{n+1}, q^{-x}, cq^{x-N}\\
aq, q^{-N}, bcq\end{array};q,q\right], \label{2.2}
\end{gather}
see also Askey and Wilson~\cite{2,3}. For $q$-analogue notations see, for example, Gasper and Rahman~\cite{9}. An idea that originated in the theory of quantum angular momentum in Physics was taken up by two mathematicians and transformed into a rich new area of research in orthogonal polynomials of a single variable. The $q$-polynomials~\eqref{2.2} and their continuous versions have been the object of great interest over the last 25 years in the f\/ield of Special Functions. They have found applications in many dif\/ferent f\/ields, including Statistical Mechanics, Quantum Group Theory, Representation Theory, Approximation Theory, and Combinatorics.

It is known in the theory of quantum angular momenta that the $6-j$ symbols are the coupling coef\/f\/icients for 3~angular momenta, and that their orthogonality property follows from the unitary nature of the coupling transformations. As a hypergeometric orthogonality this can be written in the form
\begin{gather}
\sum^N_{x=0} f_m(x) f_n(x) = \d_{m,n},\label{2.3}
\end{gather}
where the orthonormal functions $f_m(x)$ are def\/ined by
\begin{gather}
f_m(x) = (\rho (x) h_m)^{1/2} R_m(x)\label{2.4}
\end{gather}
with
\begin{gather}
\rho(x) = {\g -N+2x\o \g -N} {(\g-N, \a+1, \b+\g+1, -N)_x\o x!(\g-N-\a, -N-\b, \g +1)_x},\label{2.5}
\end{gather}
and
\begin{gather}
h_m = {(\b+1, \a+1-\g)_N\o (\a+\b+2, -\g)_N} {(\a+\b+1 + 2m)\o \a+\b+1} {(\a+\b+1, \a+1, \b+\g +1, -N)_m\o m!(\b+1, \a-\g+1, N+\a+\b +2)_m}.\label{2.6}
\end{gather}
Physicists have also given us the $9-j$ symbols, the coupling coef\/f\/icients for 4 angular momenta, including their orthogonality, which can be written
\begin{gather}
\sum_x \sum_y (2x+1) (2y+1) (2m+1) (2n+1) \!
 \left\{\!\! \begin{array}{ccc} a&b&x\\ c&d&y\\ m&n&e\end{array}\!\!\right\} \!\!\left\{ \!\!\begin{array}{ccc}
 a&b&x\\ c&d&y\\ m'&n'&e\end{array}\!\!\right\} = \d_{m,m'} \d_{n,n'},\label{2.7}
\end{gather}
where
\begin{gather}
\left\{\!\!\begin{array}{ccc} a&b&x\\ c&d&y\\ m&n&e\end{array}\!\!\right\} = \sum_k(2k+1) W(aecn;km) W(aeby;kx)
 W(by nc;kd),\label{2.8}
\\
W(aeby;kx)\nonumber\\
\quad {}= {\D(abx)\D(byk)\D(xye)\D(aek)(2a)!(a{+}b{+}e{-}y)!(a{+}b{+}e{+}y{+}1)!\o
(a{+}b{-}x)!(a{-}b{+}x)!(b{+}y{-}k)!(b{-}y{+}k)! (x{-}y{+}e)!(y{-}x{+}e)!(a{+}e{-}k)!(a{-}e{+}k)!}\nonumber\\
\qquad{}\times {} _4F_3\left[\begin{array}{c}
k-a-e, -k-a-e-1, x-a-b, -x-a-b-1\\
-2a, y-a-b-e, -y-a-b-e-1\end{array}; 1\right],\label{2.9}
\end{gather}
with the ``triangle'' function:
\begin{gather}
\D(abc) = \left\{ {(a+b-c)!(a-b+c)!(b+c-a)!\o (a+b+c+1)!}\right\}^{1/2},\label{2.10}
\end{gather}
where the implicit assumption is that $a$, $b$, $c$ satisfy the triangle inequality, and that the expressions with the factorial symbols in \eqref{2.9} and \eqref{2.10} are all nonnegative integers. The $W$-functions above are just the normalized polynomials $f_n(x)$ in dif\/ferent notation. The symmetry properties of the $6-j$ symbols enable the physicists to transform the $W$-functions in a number of dif\/ferent ways. These relations are, of course, equivalent to the Whipple formula \cite{21} for the terminating and balanced ${}_4F_3$ functions:
\begin{gather}
{}_4F_3\left[\begin{array}{c}-n, a, b, c\\ d, e, f\end{array} ;1\right]
= {(e-a, f-a)_n\o (e, f)_n} {}_4F_3\left[\begin{array}{c}-n, a, d-b, d-c,\\ d, 1+a-e-n, 1+a-f-n\end{array} ;1\right],\label{2.11}
\end{gather}

where the ``balancedness'' is indicated in the condition
\begin{gather}
a+b+c+1 = d+e+f+n.\label{2.12}
\end{gather}

Applying \eqref{2.11} several times on the $W$-function in \eqref{2.8} we reduce it to a form that is convenient for our purposes
\begin{gather}
\left\{\begin{array}{ccc} a&b&x\\ c&d&y\\ m&n&e\end{array}\right\}
 = {\D(abx) \D(cdy) \D(xye)\D (acm)\D(bdn) \D(mne)\o
(a+b-x)! (b-a+x)! (x+y-e)! (y-x+e)! (a+c-m)! (c-a+m)!}\nonumber\\
\quad{}\times {\big((2b)!\big)^2 (2c)! (a+c+n-e)! (a+b+y-e)! (b+c+y-n)!\o
(b+d-n)! (b-d+n)! (c-d+y)! (d-c+y)! (m+n-e)! (n-m+e)!}\nonumber\\
\quad{} \times {(a+c+n+e+1)! (a+b+y +e+1)\o
(c+n - b - y)!}\nonumber\\
\quad{}\times \sum_k {(2k+1) (e-a+k)! (y-b+k)! (n-c+k)! (-1)^{b+y-k}\o
(b+y-k)! (b+y + k+1)! (b-y+k)! (a+e-k)! (a-e-k)!}\nonumber\\
\quad{}\times {1\o
(c{-}n{+}k)! (a{+}e{+}k{+}1)!} {}_4F_3\!\left[\!\begin{array}{c}
k-b-y, -k-b-y-1, x-a-b, -x-a-b-1\\
e-a-b-y, -e-a-b-y-1, -2b\end{array}\!;1\right]\nonumber\\
\quad{} \times{}_4F_3\left[\begin{array}{c}
k-b-y, -k-b-y-1, n-b-d, d+n+1-b\\
-2b, n-b-c-y, n+c-b-y+1\end{array} ;1\right]\nonumber\\
\quad{} \times {}_4F_3\left[\begin{array}{c} k-c-n, -k-c-n-1, m-a-c, -m-a-c-1\\
-2c, e-a-c-n, -e-a-c-n-1\end{array};1\right].\label{2.13}
\end{gather}

For a detailed account of the $6-j$ and $9-j$ symbols see, for example, Edmonds \cite{7}. In order to identify these $9-j$ symbols as normalized orthogonal polynomials in 2 discrete variables, we replace $a+b-x$, $c+d-y$, $a+c-m$ and $b+d-n$ by $x$, $y$, $m$ and $n$, respectively, and set
\begin{gather}
a+b+c+d-e=N,\label{2.14}
\end{gather}
and assume that $N$ takes only nonnegative integer values. Now we rewrite \eqref{2.13} in a somewhat more suggestive form:
\begin{gather}
F_{m,n} (x,y;a,b,c,d)\nonumber\\
\quad{}:= \big[(2a+ 2b + 1 - 2x) (2a+ 2c + 1- 2m) (2b + 2d + 1 - 2n) (2c + 2d + 1 - 2y)\big]^{1/2}\nonumber\\
\qquad{} \times \left\{ \begin{array}{ccc}
a&b& a+b-x\\
c&d&c+d-y\\
a+c-m& b+d-n,& a+b+c+d-N\end{array}\right\}\nonumber\\
\quad{}= A_{m,n} (x,y; a, b, c,d)
\sum^{N-y}_{\ell = 0}{2\ell + 2y - 2b - 2c - 2d-1\o
2y - 2b - 2c - 2d-1} {(2y - 2b - 2c - 2d-1)_\ell \o
\ell!}\nonumber\\
\qquad{} \times {(N+y - 2a - 2b - 2c - 2d - 1, -2b, y-n-2c, y-N)_\ell (-1)^\ell\o
(2a+1+y-N, 2y-2c-2d, n+y - 2b - 2d, N+ y - 2b - 2c - 2d)_\ell}\nonumber\\
\qquad{} \times _4F_3 \left[\begin{array}{c}
-\ell, \ell + 2y - 2b - 2c - 2d-1, -x, x-2a - 2b-1\\
-2b, N+y - 2a -2b - 2c - 2d-1, y-N\end{array}; 1\right]\nonumber\\
\qquad{}\times _4F_3 \left[\begin{array}{c}
-\ell, \ell + 2y - 2b-2c-2d-1, -n, 2d+1-n\\
-2b, y-n-2c, y-n+1\end{array};1 \right]\nonumber\\
\qquad{}\times _4F_3 \left[\begin{array}{c}
n-y-\ell, n+y+\ell - 2b-2c-2d-1, -m, m-2a-2c-1\\
-2c, N+n - 2a - 2b - 2c - 2d-1, n-N\end{array};1 \right],\label{2.15}
\end{gather}
where
\begin{gather}
 A_{m,n} (x,y; a, b, c, d)= {(-y)_n\o (-N)_n} {(2b+2c+2d -N-y)!\o
(2b + 2c +2d -2y)!}\nonumber\\
 \times \Bigg\{ {N\choose m,n} {N\choose x,y} {(2a-x)! (2a-m)! (2b)! (2b)!(2c)! (2d-n)!\o
(2a+y-N)! (2a+y-N)! (2b-x)! (2b-n)! (2c-y)! (2d-y)!}\nonumber\\
 + {(2a + 2b+1 - 2x)\o
(2a + 2b+1-x)} {(2a + 2b + y - x - N)!\o
(2a + 2b-x)!}\nonumber\\
 \times {(2a + 2c +1-2m)\o
(2a + 2c+1 - m)} {(2a + 2c + n-m-N)!\o
(2a + 2c - m)!} {(2b+2d+1 - 2n)\o
(2b+ 2d+1 -n)} {(2b+2d -n-y)!\o
(2b+2d-n)!}\nonumber\\
 \times {(2c+2d +1 - 2y)\o (2c+2d +1-y)} {(2c+2d -2y)!\o
(2c+2d -y)!} {(2c+2d - 2y)! (2b+2d -n-y)!\o (2c + 2d +x-y-N)! (2b + 2d +m-n-N)!}\nonumber\\
 \times {(2a+2b + 2c+2d +1-N-n)! (2a+2b+2c + 2d +1-N-n)!\o
(2a + 2b+2c + 2d +1 - x-y-N)! (2a + 2b + 2c + 2d+1-N-m-n)!}\Bigg\}^{1/2}.\label{2.16}
\end{gather}
By repeated use of~\eqref{2.11} to transform the balanced $_4F_3$ series in \eqref{2.15} it is possible to reduce the sum over $\ell$ to a very-well-poised $_4F_3[;-1]$ series which can be summed by a standard summation formula, see Bailey \cite[4.4(3)]{5}, thereby transforming the $9-j$ symbol to a triple series. But the resulting expression is not very helpful, mainly because the
ensuing $_4F_3$ series are no longer balanced, as the $_4F_3$ series in~\eqref{2.15} are,
and hence are not easily transformable.
Besides, for the purposes of this paper reduction to a triple series is not at all useful. For some time the commonly held belief was that an orthogonal polynomial in 2~variables should be expressible as a double series, and so should be the
$9-j$ symbols.
 The results of this paper seem to indicate that it is not necessarily true. We can, in fact show, without going into detailed calculation, that even the weight function in 2 variables, or the normalization constant, need not be a compact expression. For example, when $m=0$, $n=0$, \eqref{2.15} gives us the weight function for the polynomials
\begin{gather}
w_{x,y} (a,b,c,d) = {N\choose x,y} {(2a-x)! (2a)!(2b)!(2c)!(2d)!\o
(2a+y-N)! (2a+y-N)! (2b-x)! (2c-y)!(2d-y)!}\nonumber\\
\qquad{}  \times {(2a+2b+1-2x)\o 2a+2b+1-x} {(2a+2b + y-x-N)!\o (2a + 2b-x)!} {(2a+2c -N)!\o (2a+2c)!} {2c+2d +1-2y\o
2c + 2d +1-y}\nonumber\\
\qquad{}\times {(2c + 2d -2y)!\o (2c+2d-y)!} {(2c+2d-2y)!\o (2c+2d +x-y-N)!} {(2b+2d-N)!\o (2b+2d)!}\nonumber\\
\qquad{}\times {(2a + 2b + 2c + 2d +1-N)!\o (2a+2b + 2c + 2d +1-x-y-N)!}\nonumber\\
\qquad{} \times{}_3F_2^2\left[\begin{array}{c}
x+y -N, 2a + 2b+1 + y-x-N, y-2c\\
2a+1+y-N, 2y-2c-2d\end{array};1\right].\label{2.17}
\end{gather}

Because of the self-dual character of the polynomials in \eqref{2.15} it is clear that the normalization constant is an expression similar to \eqref{2.17} with $b\leftrightarrow c$ and $x$, $y$ replaced by $m$, $n$, respectively. No matter how closely one tries to identify the 5 parameters in \eqref{2.17} with those of \eqref{2.5} it would be stretching one's imagination to think of \eqref{2.17} as a 2-dimensional extension of the simple product form that one has in~\eqref{2.5}. In general, the $_3F_2$ series in \eqref{2.17} is not summable because it is not balanced (balancedness would require the unnatural condition $2b + 2d+1 = N$). In view of this reality it is hardly surprising that the $9-j$ symbols are, at best, expressible as triple series.

There are a number of limiting cases in which the $_3F_2$ series can be summed. The case that corresponds to the Hahn polynomials \cite{13} arises when any one of the~2 parameters $a$, $c$ approaches~$\infty$ (the same, of course, is true for $b$ or $d$, but to see that one has to transform the above $_3F_2$ series f\/irst). For example, if $a \to\infty$, then, by use of Gauss' summation formula~\cite{8} one f\/inds the weight function as
\begin{gather*}
 {N\choose x,y} {(2b)!(2c)!(2d)!(2d-y)!(2b+2d -N)! (2c+2d + x-y-N)!\o
(2b-x)!(2c-y)!(2d+x-N)!^2(2b+2d)!(2c+2d-y)!}
  {2c + 2d+1\o 2c+2d+1-y},
\end{gather*}
which is a 2-dimensional extension of the weight function for the Hahn polynomials. Compare this with those, in, say~\cite{19} and~\cite{20}.

The limit case we are interested in in this paper is obtained in the limit $t\to \infty$ after setting
\begin{gather}
2a = p_1t, \qquad 2b = p_2t,\qquad 2c = p_3t,\qquad 2d = p_4t.\label{2.18}
\end{gather}
We get, as the weight function for the corresponding polynomials, the expression
\begin{gather}
 \rho_{x,y} (p_1, p_2, p_3, p_4)
 = \lim_{t\to\infty} w_{x,y} (p_1t/2, p_2t/2, p_3t/2, p_4t/2)\nonumber\\
\phantom{\rho_{x,y} (p_1, p_2, p_3, p_4)}{}  = {N\choose x,y} \eta^x_1 \eta^y_2 (1-\eta_1 - \eta_2)^{N-x-y},\label{2.19}
\end{gather}
which is a trinomial distribution, with
\begin{gather}
\eta_1  = {p_1 p_2 (p_1 + p_2 +p_3 + p_4)\o
(p_1 + p_2) (p_1 + p_3) (p_2 + p_4)},\label{2.20}\\
\eta_2  = {p_3 p_4 (p_1 + p_2 + p_3 + p_4)\o
(p_1 + p_3) (p_4 + p_2) (p_4 + p_3)}.\label{2.21}
\end{gather}
Because of the orthonormality of the $9-j$ symbols, it is guaranteed that
\begin{gather}
1-\eta_2-\eta_2 = {(p_1 p_4 - p_2p_3)^2\o
(p_1 + p_2) (p_1 + p_3) (p_4 + p_2) (p_4 + p_3)},\label{2.22}
\end{gather}
which, of course follows also from \eqref{2.20} and \eqref{2.21}. The corresponding orthonormal functions obtained from the limit of \eqref{2.15} are
\begin{gather}
R_{m,n} (x,y; p_1, p_2, p_3, p_4) \nonumber\\
 = \left\{ {N\choose x,y} {N\choose m,n}\right\}^{1/2} \Big\{ p_1^{2N-2y-x-m} p_2^{x+n} p_3^{y+m} p_4^{y-n} (p_1 + p_2)^{y-N}(p_1+p_3)^{n-N}\nonumber\\
  \times (p_2+p_4)^{N-m - 2y} (p_3 + p_4)^{N-x-2y} (p_1 + p_2 + p_3 + p_4)^{m-n+x+y}\big\}^{1/2}\nonumber\\
  \times (p_2 + p_3 + p_4)^{y-N} {(-y)_n\o (-N)_n}\nonumber\\
  \times \sum^{N-y}_{\ell = 0} {(y-N)_\ell\o \ell !} \left\{ { p_2 p_3 (p_1 + p_2 + p_3 + p_4)\o
p_1 (p_2 + p_4) (p_3 + p_4)}\right\}^\ell F\left[\begin{array}{c}
-\ell, -x\\
y-N\end{array} ; {(p_1 + p_2) (p_2 + p_3 + p_4)\o
p_2 (p_1 + p_2 + p_3 + p_4)}\right]\nonumber\\
 \times F\left[\begin{array}{c}
-\ell, -n\\ y - n+1\end{array}; - {p_4(p_2+p_3 + p_4)\o p_2p_3}\right] F\left[\begin{array}{c}
n-y-\ell, -m\\
n-N\end{array}; {(p_1 + p_3) (p_2 + p_3 + p_4)\o
p_3 (p_1 + p_2 + p_3 + p_4)}\right].\!\!\!\label{2.23}
\end{gather}
In Section \ref{section3} we will show that
\begin{gather}
R_{m,n} (x, y; p_1, p_2, p_3, p_4) \nonumber\\
\qquad{} = \left\{b_2 (x, y; N; \eta_1, \eta_2) b_2 (m, n; N; \bar\eta_1, \bar\eta_2) (1-\eta_1 - \eta_2)^{-N}\right\}^{1/2} P_{m,n} (x,y),\label{2.24}
\end{gather}
where
\begin{gather}
b_2 (x, y; N; \eta_1, \eta_2) = {N\choose x,y} \eta^x_1 \eta^y_2 (1-\eta_1 - \eta_2)^{N-x-y},\label{2.25}
\end{gather}
and
\begin{gather}
\bar\eta_1  = {p_1 p_3 (p_1 + p_2 + p_3 + p_4)\o
(p_1 + p_2) (p_1 + p_3) (p_3 + p_4)},\label{2.26}\\
\bar\eta_2  = {p_2 p_4 (p_1 + p_2 + p_3 + p_4)\o
(p_1 + p_2) (p_2 + p_4) (p_4 + p_3)},\label{2.27}
\end{gather}
(it is easily verif\/ied that $1 - \bar\eta_1 - \bar\eta_2 = 1-\eta_1 - \eta_2$),
\begin{gather}
P_{m,n} (x,y) = \sum_i\!\sum_j\!\sum_k\!\sum_\ell {(-m)_{i+j}(-n)_{k+\ell}(-x)_{i+k} (-y)_{j+\ell}\o
i!j!k!\ell! (-N)_{i+j+k+\ell}} t^iu^jv^kw^\ell,\label{2.28}
\end{gather}
with
\begin{gather}
t = {(p_1 + p_2) (p_1 + p_3)\o
p_1(p_1 + p_2 + p_3 + p_4)}, \qquad u= {(p_1 + p_3) (p_4 + p_3)\o
p_3 (p_1 + p_2 + p_3 + p_4)}, \label{2.29}\\
v  = {(p_1 + p_2) (p_2 + p_4)\o
p_2(p_1 + p_2 + p_3 + p_4)}, \qquad w = {(p_4 + p_2) (p_4 + p_3)\o
p_4 (p_1 + p_2 + p_3 + p_4)}.\nonumber
\end{gather}

One can think of the parameters $\bar\eta_1$, $\bar\eta_2$  as dual to $\eta_1$, $\eta_2$, and the polynomials $P_{m,n}(x,y)$ as self-dual. The normalization constant in $m$, $n$ is $(1-\bar\eta_1 - \bar\eta_2)^{-N}$ $b_2 (m, n; N;\bar\eta_1, \bar\eta_2)$, and that in $x, y$ is $b_2 (x,y;N; \eta_1,\eta_2)$ $(1-\eta_1-\eta_2)^{-N}$. The polynomials are very dif\/ferent from those obtainable as limits of (3.9) or (3.10) of Rahman~\cite{19}, or from the limits of the 2-dimensional case of Tratnik's~\cite{20} polynomials. Given $\eta_1$, $\eta_2$ and the trinomial distribution~\eqref{2.25} it would be almost impossible to construct the set of polynomials $P_{m,n}(x,y)$ that has 4 parameters.

After having established~\eqref{2.24} in Section~\ref{section3} we shall proceed in Section~\ref{section4} to show that the eigenfunctions $\Psi_{m,n}$ of~\eqref{1.5} are precisely the polynomials $P_{m,n}$ and the eigenvalues $\lm_{m,n}$ are certain nonlinear functions of the 4~parameters $\a_1$, $\a_2$, $\b_1$, $\b_2$ introduced in Section~\ref{section1}. These functions can be determined from symmetry considerations by requiring that for balance at equilibrium (when $\lm_{0,0} = 1$), one must have
\begin{gather}
\Psi_{0,0} (j_1, j_2)   K (i_1, i_2; j_1, j_2) = \Psi_{0,0} (i_1, i_2) K (j_1, j_2; i_1, i_2), \label{2.30}
\end{gather}
where $K$ is def\/ined in~\eqref{1.4}, and hence $\Psi_{0,0} (i_1, i_2)$ must be of the form
\begin{gather}
\Psi_{0,0} (i_1, i_2) = {N\choose i_i, i_2} \eta_1^{i_i} \eta^{i_2}_2 (1-\eta_1 - \eta_2)^{N-i_1 - i_2},\label{2.31}
\end{gather}
with
\begin{gather}
{\eta_1 (1-\a_1)\o \b_1} = {\eta_2 (1-\a_2)\o \b_2} = {1-\eta_1 - \eta_2\o 1-\beta_1 - \b_2},\label{2.32}
\end{gather}
which leads to
\begin{gather}
{\eta_1 (1-\a_1)\o \b_1} = {\eta_2 (1-\a_2)\o \b_2} = {1-\eta_1-\eta_2\o
1-\b_1 -\b_2} = D^{-1},\label{2.33}\\
D= 1+ {\a_1\b_1\o 1-\a_1} + {\a_2\b_2\o 1-\a_2}.\label{2.34}
\end{gather}

\section{A bivariate extension of Krawtchouk polynomials}\label{section3}

In order to reduce \eqref{2.23} to \eqref{2.24}--\eqref{2.29} we shall make frequent use of 3 well-known transformation formulas:
\begin{gather}
F\left[\begin{array}{c} a, b\\ c\end{array};x\right] = (1-x)^{-a} F\left[\begin{array}{c} a, c-b\\ c\end{array} ; {x\o x-1}\right]
= (1-x)^{c-a-b} F\left[\begin{array}{c} c-a, c-b\\ c\end{array} ;x\right],\label{3.1}\\
F\left[\begin{array}{c} -n, b\\ c\end{array};x\right] = {(c-b)_n\o (c)_n} F\left[\begin{array}{c} -n, b\\ 1+b-c-n\end{array};1-x\right], \qquad
n = 0, 1, 2, \ldots,\label{3.2}\\
F_1 (a; b, c;d; x,y) = (1-y)^{-a} F_1 \left(a;b, d-b-c;d; {y-x\o y-1}, {y\o y-1}\right)\nonumber\\
\phantom{F_1 (a; b, c;d; x,y)}{} = (1-x)^{-a} F_1 \left(a; d-b-c, c; d, {x\o x-1}, {x-y\o x-1}\right), \label{3.3}
\end{gather}
where
\begin{gather}
F_1 (a; b, c;d; x,y) = \sum_i\!\sum_j {(a)_{i+j}(b)_i(c)_j\o i!j!(d)_{i+j}} x^iy^j\label{3.4}
\end{gather}
is an Appell function, see, for example, Erde'lyi et al.~\cite{8}. First, by \eqref{3.1}
\begin{gather}
F\left[\begin{array}{c} -m, n-y-\ell\\
n-N\end{array}; {(p_1 + p_3) (p_2 + p_3 + p_4)\o p_3 (p_1 + p_2 + p_3 + p_4)}\right]\nonumber\\
\qquad {}= \left[- {p_1 (p_2 + p_4)\o p_3 (p_1 + p_2 + p_3 + p_4)}\right]^m F\left[\begin{array}{c} -m, \ell + y-N\\
n-N\end{array}; {(p_1 + p_3) (p_2 + p_3 + p_4)\o p_1 (p_2 + p_4)} \right],\label{3.5}
\end{gather}
and
\begin{gather}
F \left[\begin{array}{c} -x, - \ell\\
y-N\end{array}; {(p_1 + p_2) (p_2 + p_3 + p_4)\o p_2 (p_1 + p_2 + p_3 + p_4)}\right]\nonumber\\
\qquad {}= \left[ - {p_1 (p_3 + p_4)\o p_2 (p_1 + p_2 + p_3 + p_4)}\right]^x F\left[\begin{array}{c}
-x, \ell + y-N\\
y-N\end{array}; {(p_1 + p_2) (p_2 + p_3 + p_4)\o p_1 (p_3 + p_4)}\right].\label{3.6}
\end{gather}
Since
\begin{gather}
\sum_\ell {(i + y-N, j + y -N)_\ell\o \ell!(y-N)_\ell} (-\ell)_k \left\{{p_2p_3 (p_1 + p_2 + p_3 + p_4)\o
p_1 (p_2 + p_4)(p_3+p_4)}\right\}^\ell\nonumber\\
 \quad = {(i+y-N, j+y-N)_k\o (y-N)_k} \left\{ - {p_2p_3 (p_1 + p_2 + p_3 + p_4)\o p_1 (p_2 + p_4)(p_3 + p_4)}\right\}^k\nonumber\\
\qquad{} \times F\left[\begin{array}{c} i + k + y-N, j+k + y-N\\
k+y-N\end{array}; {p_2p_3(p_1+p_2 + p_3 + p_4)\o p_1 (p_2 + p_4) (p_3 + p_4)}\right]\nonumber\\
 \quad {}= {(i + y-N, j+y-N)_k\o (y-N)_k} \left\{ {p_2 p_3 (p_1 + p_2 + p_3 + p_4)\o (p_2 + p_3 + p_4) (p_2 p_3 - p_1 p_4)}\right\}^k\label{3.7}\\
\qquad{} \times \left\{ {(p_2 + p_3 + p_4) (p_1 p_4 - p_2 p_3)\o p_1 (p_2 + p_4) (p_3 + p_4)}\right\}^{N-y-i-j}
\!  F\left[\begin{array}{c}
-i, -j\\
k+y-N\end{array}\!; {p_2 p_3 (p_1 + p_2 + p_3 + p_4)\o p_1 (p_2 + p_4) (p_3 + p_4)}\right],\nonumber
\end{gather}
we have as the contribution from \eqref{3.6}
\begin{gather}
 (y-N)_k \sum_i {(-x, k+y-N)_i\o i! (y-N)_i} (-i)_\ell \left[{(p_1 + p_2) (p_2 + p_4)\o p_1 p_4 - p_2 p_3}\right]^i
 \nonumber\\
 \quad{}= {(-x)_\ell (y-N)_{k+\ell}\o (y-N)_\ell} \left[ - {(p_1 + p_2) (p_2 + p_4)\o p_1 p_4 - p_2 p_3}\right]^\ell
 \nonumber\\
 \qquad{}\times F\left[ \begin{array}{c}
\ell - x, k + \ell + y-N\\
\ell + y-N\end{array}; {(p_1 + p_2) (p_2 + p_4)\o p_1 p_4 - p_2 p_3}\right]\nonumber\\
 \quad{}= \left[- {p_2 (p_1 + p_2 + p_3 + p_4)\o p_1 p_4 - p_2 p_3}\right]^x
 \left[{(p_1 + p_2)(p_2 + p_4)\o p_2 (p_1 + p_2 + p_3 + p_4)} \right]^\ell {(-x)_\ell (y-N)_{k+\ell}\o (y-N)_\ell}
 \nonumber\\
 \qquad{}\times F\left[\begin{array}{c}
-k, \ell -x\\
\ell + y-N\end{array}; {(p_1 + p_2) (p_2 + p_4)\o p_2 (p_1 + p_2 + p_3 + p_4)}\right],\label{3.8}
\end{gather}
and that from \eqref{3.5}:
\begin{gather}
 \sum_j {(-m)_j (y-N)_{j+k}\o j!(n-N)_j} (-j)_\ell \left\{ {p_1 + p_3)(p_3 + p_4)\o p_1 p_4 - p_2 p_3}\right\}^j
 \nonumber\\
\quad{} = \left\{ - {p_3(p_1 + p_2 + p_3 + p_4)\o p_1 p_4 - p_2 p_3}\right\}^m
 \left\{{(p_1 + p_3) (p_3 + p_4)\o p_3 (p_1 + p_2 + p_3 + p_4)}\right\}^\ell {(-m)_\ell (y-N)_{k+\ell}\o (n-N)_\ell}\nonumber\\
 \qquad{}\times F\left[ \begin{array}{c}
\ell - m, n-y-k\\
\ell + n-N\end{array}; {(p_1+p_3)(p_3 + p_4)\o p_3 (p_1 + p_2 + p_3 + p_4)}\right].\label{3.9}
\end{gather}
Collecting the expressions on the right sides of \eqref{3.7}--\eqref{3.9} we f\/ind that the series part of \eqref{2.33} equals
\begin{gather}
 \left\{ {p_1 (p_3 + p_4)\o p_1p_4 - p_2 p_3}\right\}^x \left\{{(p_1 p_4 - p_2 p_3) (p_2 + p_3 + p_4)\o p_1 (p_1 + p_4) (p_3 + p_4)} \right\}^{N-y} \left\{ {p_1 (p_2 + p_4)\o p_1 p_4 - p_2 p_3}\right\}^m\nonumber\\
\qquad{} \times\sum_i\!\sum_j\!\sum_\ell {(-m)_{j+\ell} (-x)_{i+\ell}\o i!j!\ell! (y-N)_{i+\ell}(n-N)_{j+\ell}}  \left\{ {(p_1 + p_2) (p_2 + p_4) \o p_2 (p_1 + p_2 + p_3 + p_4)}\right\}^i\nonumber\\
\qquad{} \times \left\{ {(p_1 + p_3)(p_3 + p_4)\o p_3 (p_1 + p_2 + p_3 + p_4)} \right\}^j \left\{{(p_1 + p_2)(p_1 + p_3)\o p_1 (p_1 + p_2 + p_3 + p_4}\right\}^\ell S_{i,j,\ell},\label{3.10}
\end{gather}
where
\begin{gather}
 S_{i,j,\ell} =\sum_k {(-n)_k (y-N)_{k+\ell} (n-y-k)_j (-k)_i\o
 k! (y-n+1)_k} \left\{ {p_4 (p_1 + p_2 + p_3 + p_4)\o p_1 p_4 - p_2 p_3}\right\}^k\nonumber\\
\phantom{S_{i,j,\ell}}{}=(-1)^{i+j} \left\{ {p_4 (p_1 + p_2 + p_3 + p_4)\o p_1 p_4 - p_2 p_3} \right\}^i {(-n)_i (y-N)_{i+\ell}\o (y-n+1)_{i-j}}\nonumber\\
\phantom{S_{i,j,\ell}=}{}\times F\left[ \begin{array}{c} i - n, i + \ell + y-N\\
y-n+1 + i - j\end{array}; {p_4 (p_1 + p_2 + p_3 + p_4)\o p_1 p_4 - p_2 p_3}\right]\nonumber\\
\phantom{S_{i,j,\ell}}{} = (-1)^{i+j} \left\{ {p_4 (p_1 + p_2 + p_3 + p_4)\o p_1 p_4 - p_2 p_3}\right\}^\ell {(-n)_i (y-N)_{i+\ell}\o (y-n+1)_{n-j}} (N+1 - n - j- \ell)_{n-i}\nonumber\\
\phantom{S_{i,j,\ell}=}{}\times F\left[\begin{array}{c}
i-n, i + \ell + y-N\\
i + j + \ell -N\end{array}; - {(p_2 + p_4)(p_3 + p_4)\o p_1 p_4 - p_2 p_3}\right]
\nonumber\\
\phantom{S_{i,j,\ell}}{}\overset{\rm by \eqref{3.2}}{=} \left[{p_4 (p_1 + p_2 + p_3 + p_4)\o p_1 p_4 - p_2 p_3}\right]^n {(-N)_n (-n)_i (y-N)_{i+\ell} (n-N)_{j+\ell} (-y)_j\o
(-y)_n (-N)_{i+ j + \ell}}\nonumber\\
\phantom{S_{i,j,\ell}=}{}\times F \left[\begin{array}{c}
i-n, j-y\\
i+j + \ell - N\end{array}; {(p_2 + p_4)(p_3 + p_4)\o p_4 (p_1 + p_2 + p_3 + p_4)}\right].\label{3.11}
\end{gather}
Using \eqref{3.10} and \eqref{3.11} in \eqref{2.23}, and simplifying the coef\/f\/icients, we f\/inally obtain \eqref{2.24}--\eqref{2.29}. Note that, by \eqref{3.3},
\begin{gather}
P_{m,n} (x,y)  = \sum_i\!\sum_j {(-m)_{i+j} (-x)_i (-y)_j\o
i!j! (-N)_{i+j}} t^iu^j  F_1 (-n; i-x, j-y; i+j-N;v,w)\nonumber\\
\phantom{P_{m,n} (x,y)}{} = (1-v)^n \sum_i\!\sum_j {(-m)_{i+j}(-x)_i (-y)_j\o
i!j! (-N)_{i+j}} t^i u^j\nonumber\\
\phantom{P_{m,n} (x,y)=}{} \times F_1\left( -n;x + y-N, j-y; i+j -N; {v\o v-1}, {v-w\o v-1}\right)\nonumber\\
\phantom{P_{m,n} (x,y)}{}= (1-w)^n \sum_i\!\sum_j {(-m)_{i+j} (-x)_i (-y)_j\o i!j!(-N)_{i+j}} t^iu^j\nonumber\\
\phantom{P_{m,n} (x,y)=}{} \times F_1 \left(-n; i-x, x+y-N; i+j-N; {w-v\o w-1}, {w\o w-1}\right)\label{3.12}
\end{gather}
which will be very useful in the next section.

For notational simplicity let us adopt the symbol
\[
F^{(2)}_1 (a, a'; b, c; d; \lm, \mu, \nu, \rho)
\]
for the iterate of $F_1$, i.e.
\begin{gather}
F^{(2)}_1 (a, a'; b, c; d; \lm, \mu, \nu, \rho):= \sum_i\!\sum_j\!\sum_k\!\sum_\ell {(a)_{i+j} (a')_{k+\ell} (b)_{i+k}(c)_{j+\ell}\o
i!j!k!\ell! (d)_{i+j+k +\ell}} \lm^i \mu^j \nu^k \rho^\ell.\label{3.13}
\end{gather}

\section[Eigenvalues and eigenfunctions of $K (i_1, i_2; j_1, j_2)$]{Eigenvalues and eigenfunctions of $\boldsymbol{K (i_1, i_2; j_1, j_2)}$} \label{section4}

 From \eqref{1.4} and \eqref{2.25} it follows that
\begin{gather}
 K(i_1, i_2; j_1, j_2) = b_2 (i_1, i_2; N; \b_1, \b_2) (1-\a_1)^{j_1} (1-\a_2)^{j_2}\nonumber\\
\phantom{K(i_1, i_2; j_1, j_2) =}{} \times F_3 \left( -i_1, -i_2, -j_1, -j_2; -N; {\a_1\o \b_1(\a_1-1)}, {\a_2\o \b_2 (\a_2-1)}\right),\label{4.1}
\end{gather}
where the Appell function $F_3$ is def\/ined by
\begin{gather}
F_3 (a, b, a', b'; c; x, y) = \sum_r\!\sum_s {(a,a')_r (b, b')_s\o r!s! (c)_{r+s}} x^r y^s.\label{4.2}
\end{gather}
Since our objective is to show that $P_{m,n} (j_1, j_2)$ are the eigenfunctions of $K(i_1, i_2; j_1, j_2)$ for certain choices of the parameters $t$, $u$, $v$ and $w$, it suf\/f\/icies to compute the sum
\begin{gather}
Q_{m,n} (i_1, i_2) := \sum_{j_1} \sum_{j_2} b_2 (j_1, j_2; N; \eta_1, \eta_2) (1-\a_1)^{j_1} (1-\a_2)^{j_2}\label{4.3}\\
\phantom{Q_{m,n} (i_1, i_2) :=}{}  \times \! \sum_r\!\sum_s {(-i_1, -j_1)_r (-i_2, -j_2)_s\o r!s! (-N)_{r+s}} \left({\a_1\o \b_1(\a_1-1)}\right)^r \! \left({\a_2 \o \b_2 (\a_2-1)}\right)^s \!
  P_{m,n} (j_1, j_2).\nonumber
\end{gather}
We shall do the $j_1$-sum f\/irst. To facilitate the summing process we use the f\/irst of the two formulas in \eqref{3.12} to obtain
\begin{gather}
P_{m,n} (j_1, j_2)  = (1-t)^m (1-v)^n \sum_i\!\sum_j\!\sum_k\!\sum_\ell {(-m)_{i+j} (-n)_{k+\ell}\o
i!j!k!\ell!}\nonumber\\
\phantom{P_{m,n} (j_1, j_2)  =}{} \times {(j_1 + j_2 - N)_{i+k} (-j_2)_{j+\ell}\o
(-N)_{i+j + k +\ell}} \left({t\o t-1}\right)^i\! \left({t-u\o t-1}\right)^j\! \left({v\o v-1}\right)^k \! \left({v-w\o v-1}\right)^\ell.\label{4.4}
\end{gather}
Since
\begin{gather}
 \sum_{j_1} {N\choose j_1, j_2} \left( \eta_1 (1-\a_1)\right)^{j_1} \left( 1 - \eta_1 - \eta_2\right)^{N-j_1 - j_2} (-j_1)_r (j_1 + j_2 - N)_{i+k}\nonumber\\
\qquad {}= {N\choose j_2}  \left( \eta_1 (1-\a_1)\right)^r (1-\eta_1 - \eta_2)^{i+k} (1-\a_1 \eta_1 - \eta_2)^{N-j_2 - r - i - k} (j_2 - N)_{r+i+k},\!\!\!\label{4.5}
\end{gather}
the r.h.s.\ of \eqref{4.3} can be written as
\begin{gather}
 (1-t)^m(1-v)^n \sum_{j_2} {N\choose j_2} (\eta_2 (1-\a_2))^{j_2} (1-\a_1\eta_1 - \eta_2)^{N-j_2}\nonumber\\
\qquad {}\times \sum_r\!\sum_s {(-i_1, j_2 - N)_r (-i_1, - j_2)_s\o
r!s! (-N)_{r+s}} \left( - {\a_1\eta_1\o \b_1 (1-\a_1 \eta_1 - \eta_2)}\right)^r\nonumber\\
\qquad {}\times \left({\a_2\o \b_2 (\a_2-1)}\right)^s F^{(2)}_1 \Bigg( -m, -n; r+j_2 -N, - j_2; -N; {t(1-\eta_1-\eta_2)\o (t-1)(1-\a_1\eta_1 - \eta_2)},\nonumber\\
\qquad \quad \left({t-u\o t-1}\right), \left({v(1-\eta_1 - \eta_2)\o (v-1)(1-\a_1 \eta_1 - \eta_2)}\right), \left({v-w\o v-1}\right)\Bigg),\label{4.6}
\end{gather}
by using the transformation formulas
\begin{gather}
 F^{(2)}_1 (a, a'; b, c; d; \lm, \mu, \nu, \rho)\nonumber\\
\qquad{}= (1-\lm)^{-a} (1-\nu)^{-a'} F^{(2)}_1 \left(a, a'; d-b-c, c;d; {\lm\o \lm-1}, {\lm -\mu \o \lm-1}, {\nu\o \nu-1}, {\nu-\rho\o \nu-1}\right)\label{4.7}\\
\qquad{}= (1-\mu)^{-a} (1-\rho)^{-a'} F^{(2)}_1 \left(a, a'; b, d-b-c;d; {\mu-\lm\o \mu-1}, {\mu\o \mu-1}, {\rho-\nu\o \rho-1}, {\rho\o \rho-1}\right),\label{4.8}
\end{gather}
which are direct consequences of \eqref{3.3}. By~\eqref{4.8}, the $F^{(2)}_1$ series in~\eqref{4.6} transforms to
\begin{gather}
\left({1-u\o 1-t}\right)^m \left({1-w\o 1-v}\right)^n F^{(2)}_1 (-m, -n; r+ j_2 - N, -r; -N; t', u', v', w'),\label{4.9}
\end{gather}
where
\begin{gather}
t'  = {(t-u) (1-\a_1 \eta_1 - \eta_2) - t (1-\eta_1 - \eta_2)\o (1-u) (1-\a_1\eta_1 - \eta_2)} = {t\b_1 - u (1-\b_2)\o (1-u) (1-\b_2)},\label{4.10}\\
u'  = {t-u\o 1-u}, \label{4.11}\\
v'  = {(v-w)(1-\a_1 \eta_1 - \eta_2) - v(1 -\eta_1 - \eta_2)\o (1-w) (1-\a_1 \eta_1 - \eta_2)} = {v\b_1 - w(1-\b_2)\o (1-w)(1-\b_2)},\label{4.12}\\
w'  = {v-w\o 1-w}, \label{4.13}
\end{gather}
where the expressions for $t'$, $v'$ in terms of $\b_1$, $\b_2$ are obtained from \eqref{2.32} and \eqref{2.33}. With~\eqref{4.9} we may now carry out the sum over $j_2$. We have
\begin{gather}
 \sum_{j_2} {N\choose j_2} \left(\eta_2 (1-\a_2)\right)^{j_2} (1-\a_1 \eta_1 - \eta_2)^{N-j_2} (-j_2)_s (j_2 -N)_{r+i+k}\nonumber\\
\qquad{}= \left({\eta_1 (1-\a_2)\o 1-\a_1 \eta_1 - \a_2 \eta_2}\right)^s \left({1-\a_1 \eta_1 - \eta_2\o 1-\a_1 \eta_1 - \a_2 \eta_2}\right)^{r+i +k} (1-\a_1 \eta_1 - \a_2 \eta_2)^N(-N)_{r+s+i+k}\nonumber\\
\qquad{}= \b^s_2 (1-\b_2)^{r+i+k} D^{-N} (-N)_{r+s+i+k}.\label{4.14}
\end{gather}
So the expression in \eqref{4.6} can be written as
\begin{gather}
 (1-u)^m(1-w)^n (1-\a_1 \eta_1 - \a_2 \eta_2)^N \sum_r\!\sum_s {(-i_1)_r (-i_2)_s \o r!s!} \left({\a_1\o \a_1 -1}\right)^r \left({\a_2 \o \a_2 -1}\right)^s\label{4.15}\\
{}\times F^{(2)}_1 \left( -m, -n; r+s-N, -r; -N; {t\b_1 -u(1-\b_2)\o (1-u)(1-\b_2)}, {t-u\o 1-u}, {v\b_1 - w(1-\b_2)\o (1-w) (1-\b_2)}, {v-w\o 1-w}\right).\nonumber
\end{gather}
We are just one transformation away from the form where we can do the $r$ and $s$ summations. First, we interchange the f\/irst and third parameters, then the second and the fourth, so that the parameters $r+s-N$ and $-r$ are also interchanged (although this step is not necessary), followed by an application of \eqref{4.8}. The end result is that \eqref{4.15} transform to
\begin{gather}
 (1-\a_1\eta_1 - \a_2 \eta_2)^N \left\{{1-\a_1\eta_1 - \a_2\eta_2 - \eta_1t(1-\a_1) - \eta_2 u(1-\a_2)\o
1-\a_1\eta_1 - \a_2 \eta_2}\right\}^m\nonumber\\
\qquad{}\times\left\{{1-\a_1\eta_1 - \a_2\eta_2 - \eta_1 v(1-\a_1) - \eta_2 w (1-\a_2)\o
1-\a_1\eta_1 - \a_2\eta_2}\right\}^n\label{4.16}\\
 \qquad{}\times \sum_r\!\sum_s {(-i_1)_r (-i_1)_s\o r!s!} \left({\a_1\o \a_1-1}\right)^r\left({\a_2\o \a_2-1}\right)^s F_1^{(2)} (-m, -n; -r, -s; -N; \lm, \mu, \nu, \rho),\nonumber
\end{gather}
with
\begin{gather}
\lm  = {(1-\b_1) t - \b_2 u\o 1-\b_1 t - \b_2 u},\label{4.17}\\
\mu  = {u(1-\b_2) - \b_1t\o 1-\b_1t - \b_2 u}, \label{4.18}\\
\nu  = {\eta_2 (\a_2-1) (v-w) - v (1-\eta_1 - \eta_2)\o
v \eta_1 (1-\a_1) + w \eta_2 (1-\a_2) - (1-\a_1 \eta_1 - \a_2 \eta_2)},\label{4.19}\\
\rho  = {(1-\a_1\eta_1- \a_2 \eta_2) (v-w) - v (1-\eta_1 - \eta_2)\o
v \eta_1 (1-\a_1) + w \eta_2 (1-\a_2) - (1-\a_1\eta_1 - \a_2\eta_2)}.\label{4.20}
\end{gather}

Now,
\[
\sum_r {(-i_1)_r\o r!} \left({\a_1 \o \a_1 -1}\right)^r (-r)_{i+k} = \a_1^{i+k} (1-\a_1)^{-i_1}(-i_1)_{i+k},
\]
and
\[
\sum_s {(-i_2)_s\o s!} \left({\a_2\o \a_2 -1}\right)^s  (-s)_{j+\ell} = \a_2^{j+\ell} (1-\a_2)^{-i_2} (-i_2)_{j+\ell},
\]
so that the double sum in \eqref{4.16} reduces to
\begin{gather}
(1-\a_1)^{-i_1}(1-\a_2)^{-i_2} F^{(2)}_1 (-m, -n; -i_1, -i_2; -N; \a_1 \lm, \a_2\mu, \a_1\nu, \a_2 \rho).\label{4.21}
\end{gather}

Note that
\begin{gather}
 b_2 (i_1,i_2; N; \b_1, \b_2) (1-\a_1)^{-i_1} (1-\a_2)^{-i_2} (1-\a_1 \eta_1 - \a_2\eta_2)^N\nonumber\\
\qquad{} = {N\choose i_1, i_2} \b_1^{i_1} \b_2^{i_2} (1-\b_1 - \b_2)^{N-i_1-i_2} \left({\eta_1\o
\b_1 (1-\a_1\eta - \a_2\eta_2)}\right)^{i_1} \nonumber\\
\qquad\quad{}\times \left({\eta_2\o \b_2 (1-\a_1\eta_1 - \a_2\eta_2)}\right)^{i_2}
  (1-\a_1\eta_1 - \a_2 \eta_2)^N\nonumber\\
\qquad{} = {N\choose i_1, i_2} \eta^{i_1}_1 \eta^{i_2}_2 (1-\eta_1 - \eta_2)^{N-i_1-i_2}
  = b_2 (i_1, i_2; N; \eta_1, \eta_2) \quad \mbox{by \eqref{2.32}--\eqref{2.34}}.\label{4.22}
\end{gather}
So the polynomial $P_{m,n} (i_1, i_2)$ is indeed an eigenfunction of $K (i_1, i_2; j_1, j_2)$ with the corres\-ponding eigenvalue
\begin{gather}
\lm_{m,n}  = \left\{{1-\a_1 \eta_1 - \a_2 \eta_2 - \eta_1 t (1-\a_1) - \eta_2 u (1-\a_2)\o
1-\a_1 \eta_1 - \a_2\eta_2}\right\}^m\nonumber\\
\phantom{\lm_{m,n}  =}{}
 \times \left\{{1-\a_1 \eta_1 - \a_2 \eta_2 - \eta_1 v (1-\a_1) - \eta_2 w (1-\a_2)\o
1-\a_1 \eta_1 - \a_2 \eta_2}\right\}^n\nonumber\\
\phantom{\lm_{m,n}}{} = (1-\b_1t - \b_2 u)^m (1-\b_1 v - \b_2 w)^n,\label{4.23}
\end{gather}
provided we can f\/ind solutions for the parameters $t$, $u$, $v$, $w$ in terms of $\a_1$, $\a_2$, $\b_1$, $\b_2$ such that
\begin{gather}
t = \lm \a_1, \qquad u = \mu \a_2,\qquad v = \nu \a_1, \qquad w =\rho \a_2.\label{4.24}
\end{gather}
This is, of course, elementary algebra, which will be carried out in the next two sections.

\section[Eigenvalues and eigenfunctions: the nondegenerate case $\alpha_1 \ne \alpha_2$]{Eigenvalues and eigenfunctions:\\ the nondegenerate case $\boldsymbol{\a_1 \ne \a_2}$}\label{section5}

The relations between the parameters, i.e., \eqref{4.24}, can be expressed in the form
\begin{gather}
t  = \a_1 {t(1-\b_1) - \b_2 u\o 1-\b_1t - \b_2 u},\label{5.1}\\
u  = \a_2 {u(1-\b_1) - \b_1 t\o 1-\b_1 t - \b_2 u},\label{5.2}\\
v  = \a_1 {v(1-\b_1) - \b_2u\o 1-\b_1 v - \b_2 w}, \label{5.3}\\
w  = \a_2 {w(1-\b_2) - \b_1 v \o 1-\b_1 v - \b_2 w}.\label{5.4}
\end{gather}
From \eqref{5.1}
\begin{gather}
1-\b_1t - \b_2 u = {\a_1 (1-t)\o \a_1 -t}.\label{5.5}
\end{gather}
which, on substitution in \eqref{5.2} gives
\begin{gather}
\b_2 u = - \b_1 t + {t(1-\a_1)\o t-\a_1}.\label{5.6}
\end{gather}
Combination of \eqref{5.6} and \eqref{5.1} or \eqref{5.2} gives the quadratic relation for $t$:
\begin{gather}
\b_1 (\a_1 - \a_2) (t-\a_1)^2 - (1-\a_1) (\a_1 - \a_2 +\a_1\b_1 + \a_2\b_2) (t-\a_1) + \a_1 (1-\a_1)^2 = 0.\label{5.7}
\end{gather}
If $\a_1 > \a_2$, then both roots are real and positive with the discriminant $\D$ given by
\begin{gather}
\D  = (\a_1 - \a_2 + \a_1 \b_1 + \a_2 \b_2)^2 - 4\a_1 \b_1 (\a_1 - \a_2)\nonumber\\
\phantom{\D}{} = (\a_1 - \a_2 + \a_2 \b_2 - \a_1 \b_1)^2 + 4\a_1 \a_2 \b_1 \b_2,\label{5.8}
\end{gather}
which is $>0$ since the parameters, being probabilities, are necessarily in $(0,1)$. The roots are
\begin{gather}
t - \a_1 = \big( \a_1 - \a_2 + \a_1 \b_1 + \a_2 \b_2 \pm \D^{1/2}\big) {(1-\a_1)\o 2 (\a_1 - \a_2)\b_1}.\label{5.9}
\end{gather}
From \eqref{5.1} and \eqref{5.3} it is clear that $v$ satisf\/ies the same equation as \eqref{5.7}, so we may take $t-\a_1$ and $v -\a_1$, having one of the signs indicated in \eqref{5.9} (it is immaterial which sign we assign to each). However, if $\a_1 < \a_2$, the equation is
\begin{gather}
u-\a_2 = \big(\a_2 -\a_1 +\a_1 \b_1 + \a_2\b_2 \pm\D^{1/2}\big) {(1-\a_2)\o 2\b_2 (\a_2 - \a_1)}.\label{5.10}
\end{gather}
Since $K$ is a transition probability it is necessarily positive and less than~1, and hence $0 < \lm_{m,n} < 1$. This implies that $1 - \b_1 t - \b_2 u$ and $1-\b_1 v - \b_2 w$ must both be in $(0,1)$. From \eqref{5.1}--\eqref{5.4} it is clear that
\begin{gather}
\lm_{m,n} = \left(\a_1 \left({1-t\o \a_1 -t}\right)\right)^m \left( \a_2 \left({1-w\o \a_2 -w}\right)\right)^n.\label{5.11}
\end{gather}

Obviously we have to choose values of $t$ and $w$ such that
\begin{gather}
{\rm (i)}\quad {\rm either}\quad  t > 1\quad {\rm or}\quad t < \a_1,\qquad
{\rm (ii)}\quad {\rm either}\quad   w > 1 \quad {\rm or}\quad w < \a_2. \label{5.12}
\end{gather}

A straightforward calculation, however, shows that in both cases the end-result is the same, i.e.
\begin{gather}
\a_1 {1-t\o \a_1 - t} = \a_2 {1-w\o \a_2 - w} = {1\o 2} \big\{ \a_1 (1-\b_1) + \a_2 (1-\b_2) + \D^{1/2}\big\},\label{5.13}
\end{gather}
irrespective of whether or not $\a_1 - \a_2$ is $+ve$ or $-ve$. So the eigenvalues are
\begin{gather}
\lm_{m,n} = \left\{ {\a_1(1-\b_1) + \a_2(1-\b_2) + \D^{1/2}\o 2}\right\}^{m+n}.\label{5.14}
\end{gather}

\section[The degenerate case $\alpha_1 = \alpha_2$]{The degenerate case $\boldsymbol{\a_1 = \a_2}$}\label{section6}

When $\a_1 = \a_2 = \a$, say, we subtract \eqref{5.2} from \eqref{5.1} to get
\begin{gather}
t-u = \a {t-u\o 1-\b_1 t - \b_2 u}.\label{6.1}
\end{gather}
So, either $t = u$ or $1-\b_1 t - \b_2 u = \a$. Since we must assume that $0 < \a < 1$, the second alternative is impossible as it would require $\a=1$, which can be seen from either \eqref{5.1} or \eqref{5.2}. So we must conclude that
\begin{gather}
t = u,\label{6.2}
\end{gather}
in which case it follows that $\D^{1/2} = \a(\b_1 + \b_2)$, so from \eqref{5.14}, we have
\begin{gather}
\lm_{m,n} = \a^{m+n}.\label{6.3}
\end{gather}
Similarly it follows that $v = w$, and that, ultimately
\begin{gather}
 t= u = v = w.\label{6.4}
\end{gather}

The eigenfunction in this degenerate case reduces to
\begin{gather}
 F^{(2)}_1 (-m, -n; -i_1, -i_2; -N; t, t, t, t)\nonumber\\
\qquad = \sum_r\!\sum_s {(-m)_{r+s}(-i_i)_r (-i_2)_s\o r!s! (-N)_{r+s}} t^{r+s}
  F_1 (-n; r-i_1, s-i_2; r+s - N; t,t)\nonumber\\
\qquad {} = \sum_r\!\sum_s {(-m)_{r+s} (-i_1)_r (-i_2)_s \o r!s! (-N)_{r+s}} t^{r+s}  {} _2F_1(-n, r+s -i_1 - i_2; r+s - N;t)\nonumber\\
\qquad \overset{\rm \mbox{\scriptsize by~\cite[9.5(1)]{4}}}{=}
  (1-t)^n \sum_r\!\sum_s {(-m)_{r+s}(-i_1)_r (-i_2)_s\o r!s! (-N)_{r+s}} t^{r+s}\nonumber\\
\qquad\quad {}\times {}_2F_1 \left(-n, i_1 + i_2 -N; r+s -N; {t\o t-1}\right)\nonumber\\
\qquad \overset{\rm by~\eqref{3.1}}{=}
 (1-t)^n \sum_k {(-n, i_1 + i_2 - N)_k \o k!(-N)_k} \left({t\o t-1}\right)^k F_1 (-m; -i_1, -i_2; k-N;t,t)\nonumber\\
\qquad{} = (1-t)^{m+n} F_1 \left( i_1 + i_2 - N; -m, -n; -N; {t\o t-1}, {t\o t-1}\right)\nonumber\\
\qquad{}= (1-t)^{m+n} {}_2F_1\left( -m-n, i_1 + i_2 - N; -N; {t\o t-1}\right)\nonumber\\
\qquad{}={}_2F_1 (-m-n, -i_1 - i_2; -N; t).\label{6.5}
\end{gather}
It also follows from \eqref{5.7} that
\begin{gather}
t = {1-\a (1-\b_1 - \b_2)\o \b_1 + \b_2}.\label{6.6}
\end{gather}
So, in this special case, the eigenfunction is essentially a single-variable Krawtchouk polynomial of degree $m + n$ in $i_1 +  i_2$.

\section{Concluding remarks and acknowledgements}\label{section7}

It is nearly 3 years since the f\/irst draft of this paper was prepared, then sent away for private circulation to professional friends and colleagues. Since then it was brought to our attention that a number of publications exist in the literature of both orthogonal polynomials and probability-statistics, that are closely related to what we have done in this paper. The earliest among them was, as far as we know, the 1971 paper of R.C.~Grif\/f\/iths~\cite{10}, see also \cite{11} and~\cite{12}, where he considers a (persumably more general) class of transition density expansions of the so-called Lancaster type. He used probability generating functions to characterize bivariable distributions with identical multinominal marginals, with the transition density having orthogonal polynomials as eigenfunctions, quite akin to what we have attempted to do here. One might argue, as one of the referees of our paper has pointed out, that Grif\/f\/iths' paper says more about the probabilistic nature of the model than ours do. However, our principal motivation in delving into the Quantum Angular Momentum literature is to get a handle on the problem of how to f\/it the four probability parameters into a trinomial-distribution-based cumulative Bernoulli model that has only two independent parameters. $9-j$ symbols of Angular Momentum theory provided us with a 4-parameter representation of the two probabilities in the trinomial distribution.

Fortunately, one of us (MR) had the good fortune of meeting Dr.~Grif\/f\/iths in an Orthogonal Polynomial meeting in France in 2007, and had the benef\/it of a fruitful discussion on what we had done in our paper, and what he had done much earlier. The authors gratefully acknowledge the help he provided us with reprints of his papers. We are also grateful to the f\/irst referee for pointing out the importance of discussing Dr.~Grif\/f\/iths' work in this paper, which we had intended to do in a subsequent publication.

There was yet another eye-opening experience for MR when he met Dr.~Zhedanov of Donetsk Institute for Physics and Technology in Ukraine and talked about the present paper. It turned out that Dr. Zhedanov~\cite{24} also found a very similar 2-variable Krawtchouk polynomial by considering the oscillator algebra of the $9-j$ symbols. However, his polynomials are not exactly the same as ours, but a limiting case of. It is obvious that he would have found the same polynomials as we have, had he chosen to work with the full $SU(2)$ algebra of the $9-j$ symbols.

In the latest SIDE8 meeting in Montreal, June'08, Professor M. Noumi pointed out to MR that a multidimensional version of our 2-variable Krawtchouk polynomial was found by Aomoto and Gelfand~\cite{1}, and later by Mizukawa~\cite{17}, who gave a zonal spherical functions proof of the orthogonality of the polynomials. We owe our gratitude to Dr.~Noumi as well.

\pdfbookmark[1]{References}{ref}
\LastPageEnding

\begin{thebibliography}{99}

\footnotesize\itemsep=0pt

\bibitem{1} Aomoto K., Kita M.,  Theory of hypergeometric functions, Springer, Tokyo, 1994 (in Japanese).

\bibitem{2} Askey R.,   Wilson  J.A., A set of orthogonal polynomials that generalize the Racah coef\/f\/icients or $6-j$ symbols, {\it SIAM J. Math. Anal.}  {\bf 10} (1979), 1008--1016.

\bibitem{3} Askey R., Wilson J.A., A set of hypergeometric orthogonal polynomials, {\it SIAM J. Math. Anal.} {\bf 13} (1982), 651--655. 

\bibitem{4} Askey R., Wilson J.A., Some basic hypergeometric polynomials that generalize Jacobi polynomials, {\it Mem. Amer. Math. Soc.} {\bf 319} (1985), 1--55.

\bibitem{5} Bailey W.N.,  Generalized hypergeometric series, Cambridge University Press, Cambridge,  1935 (reprinted by Stechert-Hafner, New York, 1964).

\bibitem{6} Cooper R.D., Hoare M.R.,  Rahman M., Stochastic processes and special functions: on the probabilistic origin of some positive kernels associated with classical orthogonal polynomials, {\it J. Math. Anal. Appl.}  {\bf 61} (1977), 262--291.

\bibitem{7}Edmonds A.R., Angular momentum in quantum mechanics, 2nd ed., Princeton University Press, Princeton, New Jersey, 1960.

\bibitem{8}Erd\'elyi  et. al., Higher transcendental functions, Vol.~I, McGraw-Hill, New York, 1953.

\bibitem{9} Gasper G.,   Rahman  M., Basic hypergeometric series, 2nd ed.,
{\it Encyclopedia of Mathematics and Its Applications}, Vol.~96, Cambridge University Press, Cambridge,    2004.

\bibitem{10} Grif\/f\/iths R.C., Orthogonal polynomials on the multinomial distribution, {\it Austral. J. Statist.} {\bf 13} (1971), 27--35,   Corregenda, {\it Austral. J. Statist.} {\bf 14} (1972), 270.

\bibitem{11} Grif\/f\/iths R.C., Orthogonal polynomials on the negative multinomial distribution, {\it J. Multivariate Anal.} {\bf 5} (1975), 271--277.

\bibitem{12} Grif\/f\/iths R.C., Orthogonal polynomials on the multinomial, Notes: Version 3.0, 04/09/2006, unpublished (Private communication).

\bibitem{13}  Hahn W., \"Uber Orthogonalpolynome, die $q$-Dif\/ferenzengleichungen gen\"ugen, {\it Math. Nachr.} {\bf 2} (1949), 4--34.

\bibitem{14}Hoare M.R., Rahman  M., Cumulative Bernoulli trials and Krawtchouk processes, {\it Stochastic Process. Appl.} {\bf 16} (1983), 113--139.

\bibitem{15} Hoare M.R.,   Rahman M., Cumulative hypergeometric processes: a statistical role for the $_nF_{n-1}$ functions, {\it J. Math. Anal. Appl.}  {\bf 135}  (1988), 615--626. 

\bibitem{16} Hoare M.R.,   Rahman M., Distributive processes in discrete systems, {\it Phys.~A}  {\bf 97} (1979), 1--41. 

\bibitem{17} Mizukawa H., Zonal spherical functions on the complex ref\/lection groups and $(n+1, m+1)$-hypergeometric functions, {\it Adv. Math.} {\bf 184} (2004), 1--17.

\bibitem{18} Racah G., Theory of complex spectra. II, {\it Phys. Rev.} {\bf 62} (1942), 438--462.

\bibitem{19} Rahman M., Discrete orthogonal systems corresponding to Dirichlet distribution, {\it Utilitas Math.}  {\bf 20} (1981), 261--272.

\bibitem{20} Tratnik M.V., Some multivariable orthogonal polynomials of the Askey tableau-discrete families, {\it J. Math. Phys.} {\bf 32} (1991), 2337--2342.

\bibitem{21} Whipple F.J.W., Well-poised series and other generalized hypergeometric series, {\it Proc. Lond. Math. Soc. (2)} {\bf 25} (1926), 525--544.

\bibitem{22} Wilson J.A., Hypergeometric series, recurrence relations and some new orthogonal functions, Thesis, Univ. of Wisconsin, Madison, 1978.

\bibitem{23} Wilson J.A., Some hypergeometric orthogonal polynomials, {\it SIAM J. Math. Anal.} {\bf 11} (1980), 690--701. 

\bibitem{24} Zhedanov A., $9j$-symbols of the oscillator algebra and Krawtchouk polynomials in two variables, {\it J.~Phys.~A: Math. Gen.} {\bf 30} (1997), 8337--8353.

\end{thebibliography}
\end{document}